\documentclass[reqno, letterpaper, oneside]{article}
\usepackage[utf8]{inputenc}
\usepackage[T1]{fontenc}
\usepackage{comment}
\usepackage{authblk}
\usepackage{relsize}
\usepackage{changepage}
\usepackage{mathtools}
\usepackage[english]{babel}
\usepackage{amsmath}
\usepackage{amsthm}
\usepackage{amsfonts}
\usepackage{amssymb}
\usepackage{graphicx}
\usepackage[normalem]{ulem}
\usepackage[usenames,dvipsnames]{color}
\theoremstyle{plain}
\newtheorem{theorem}{Theorem}[section]
\newtheorem{corollary}[theorem]{Corollary}

\newtheorem{remark}[theorem]{Remark}
\newtheorem{lemma}[theorem]{Lemma}

\newtheorem*{claim}{Claim}
\newtheorem{observation}[theorem]{Observation}

\makeatletter
\newcommand{\vast}{\bBigg@{4}}
\newcommand{\Vast}{\bBigg@{5}}
\makeatother
\definecolor{bulgarianrose}{rgb}{0.28, 0.02, 0.03}
\definecolor{gray}{rgb}{0.5, 0.5, 0.5}

\theoremstyle{definition}

\theoremstyle{remark}
\newtheorem*{fact*}{Fact}
\newtheorem*{question*}{Question}
\usepackage[left=2cm,right=2cm,top=2cm,bottom=2cm]{geometry}
\usepackage{amssymb}
\usepackage{hyperref}
\makeatletter
\def\namedlabel#1#2{\begingroup
    #2%
    \def\@currentlabel{#2}%
    \phantomsection\label{#1}\endgroup
}
\usepackage[normalem]{ulem}
\usepackage{dsfont}
\usepackage{accents}
\usepackage{verbatim}
\usepackage{ulem}
\usepackage{pgf,tikz,pgfplots}
\pgfplotsset{compat=1.16}
\usepackage[all]{xy} 

\usepackage{stackengine}
\stackMath
\newcommand\tsup[2][2]{%
 \def\useanchorwidth{T}%
  \ifnum#1>1%
    \stackon[-.5pt]{\tsup[\numexpr#1-1\relax]{#2}}{\scriptscriptstyle\sim}%
  \else%
    \stackon[.5pt]{#2}{\scriptscriptstyle\sim}%
  \fi%
}

\usepackage[toc,page]{appendix}

\newcommand{\boxi}{\text{box}}

\newcommand{\marco}{\textcolor{blue}}

\usepackage{subfig}

\usepackage{enumerate}

\title{\scshape
  On the boxicity of Kneser graphs and complements of line graphs}

\author[1]{Marco Caoduro}
\author[2]{Lyuben Lichev}
\affil[1]{Univ. Grenoble Alpes, Laboratoire G-SCOP, Grenoble, France}
\affil[2]{Ecole Normale Sup\'erieure de Lyon, Lyon, France}

\begin{document}

\maketitle
 
\begin{abstract}
An axis-parallel $d$-dimensional box is a cartesian product $I_1\times I_2\times \dots \times I_b$ where
$I_i$ is a closed sub-interval of the real line. For a graph $G = (V,E)$, the \emph{boxicity of $G$}, denoted by $\boxi(G)$, is the minimum dimension $d$ such that $G$ is the intersection graph of a family $(B_v)_{v\in V}$ of $d$-dimensional boxes in $\mathbb R^d$.

Let $k$ and $n$ be two positive integers such that $n\geq 2k+1$. The \emph{Kneser graph} $Kn(k,n)$ is the graph with vertex set given by all subsets of $\{1,2,\dots,n\}$ of size $k$ where two vertices are adjacent if their corresponding $k$-sets are disjoint. In this note we derive a general upper bound for $\boxi(Kn(k,n))$, and a lower bound in the case $n\ge 2k^3-2k^2+1$, which matches the upper bound up to an additive factor of $\Theta(k^2)$. Our second contribution is to provide upper and lower bounds for the boxicity of the complement of the line graph of any graph $G$, and as a corollary we derive that $\boxi(Kn(2,n))\in \{n-3, n-2\}$ for every $n\ge 5$. 
\end{abstract}

\hspace{1em}Keywords: Boxicity, Kneser Graphs, Line Graphs, Graph Theory

\hspace{1em}MSC Class: 05C62

\section{Introduction}

An axis-parallel $d$-dimensional box is a cartesian product $I_1\times I_2\times \dots \times I_b$ where
$I_i$ is a closed sub-interval of the real line. For a graph $G = (V,E)$, the \emph{boxicity of $G$}, denoted by $\boxi(G)$, is the minimum dimension $d$ such that $G$ is the intersection graph of a family $(B_v)_{v\in V}$ of $d$-dimensional boxes in $\mathbb R^d$. Boxicity has been introduced by Roberts~\cite{Rob} in 1969 and has been extensively studied since then, see for example \cite{ABC, CFS, CS, Esp, EJ, SW}.

Let $k$ and $n$ be two positive integers such that $n\geq 2k+1$. The \emph{Kneser graph} $Kn(k,n)$ is the graph with vertex set given by all subsets of $[n] := \{1,2,\dots,n\}$ of size $k$ where two vertices are adjacent if their corresponding $k$-sets are disjoint. The notion of Kneser graph was born in 1955 in a paper of Kneser~\cite{1955_Kneser}, where he conjectured that the chromatic number $\chi(Kn(k,n))$ is equal to $n - 2k + 2$. In 1978 Lovasz \cite{1978_Lovasz} settled this conjecture with a brilliant topological proof. Since then, several papers have focused on the properties of this family, see for example \cite{2000_Chen, 1986_Frankl, 2013_Harvey}.

In this note we are interested in deriving bounds on the boxicity of Kneser graphs. In particular, we prove the following theorems.

\begin{theorem} \label{main thm}
Fix two positive integers $k,n$ with $n\ge 2k+1$. The boxicity of the Kneser graph $Kn(k,n)$ is at most $n-2$. Moreover, if $n\geq 2k^3 - 2k^2 + 1$, then $\boxi(Kn(k,n))\ge n - \dfrac{13k^2-11k+16}{2}$.
\end{theorem}

In general, less precise lower bounds can be obtained without the assumption $n\geq 2k^3 - 2k^2 + 1$ by exploiting the relationship between boxicity and poset dimension proved by \cite{ABC}. We will have a quick glance at the main technique at the end of Section \ref{sec lower bounds} and we invite the reader to consult the papers~\cite{ABC, Fur, Kie1} for further details.

The second main part of this paper deals with the boxicity of complements of line graphs and in particular the boxicity of the Kneser graph $Kn(2,n)$. The \emph{line graph} $L(G)$ of a graph $G$ has vertex set $E(G)$ and edge set $\{(uv, wv):\hspace{0.2em} uv, wv\in E(G)\}$. Denote by $G^c$ the complement graph of a graph $G$, by $\delta(G)$ the minimum degree of $G$, and by $\Delta(G)$ the maximum degree of $G$.

The aim of the next theorem is twofold: first, it gives a sharper lower bound on $\boxi(Kn(2,n))$ than Theorem~\ref{main thm}, and second, it generalises this lower bound to complements of line graphs, which are realised as induced subgraphs of $Kn(2,n)$. 

\begin{theorem}\label{thm complements line graphs}
Let $G$ be any graph on $n$ vertices of maximum degree $\Delta = \Delta(G) \ge 3$ and let $H = L(G)$ denote its line graph. Then, the boxicity of the complement of $H$ is at most $n-2$. Moreover,
\begin{itemize}
    \item  $\boxi(H^c) \geq \dfrac{|E(H)|}{12}$, if $\Delta = 3$;
    \item  $\boxi(H^c) \geq \dfrac{|E(H)|}{16}$, if $\Delta = 4$;
    \item  $\boxi(H^c) \geq \dfrac{2|E(H)|}{\Delta^2+3\Delta}$, if $\Delta \geq 5$.
\end{itemize}
\end{theorem}

\begin{remark}
If $\Delta(G) = 2$, then both $G$ and $H$ are unions of disjoint paths and cycles. For this particular kind of graphs, Corollary 3.3 in \cite{1983_Cozzens} together with Lemma 3 in \cite{1979_Trotter} imply that $box(H^c) = \sum_{i=1}^k \left\lceil \frac{|E(H_i)|}{3} \right\rceil$, where $(H_i)_{i \in [k]}$ are the connected components of $H$ .
\end{remark}

\begin{corollary} \label{Box_Kneser}
For every $n\ge 5$, the boxicity of the Kneser graph $Kn(2,n)$ is either $n-3$ or $n-2$.
\end{corollary}
\begin{proof} [Proof of Corollary~\ref{Box_Kneser} assuming Theorem~\ref{thm complements line graphs}]
Note that $Kn(2,n)$ is the complement of the line graph of the complete graph $K_n$. Then, Theorem~\ref{thm complements line graphs} applied for $G = K_n$ with $\Delta = n-1$ shows the upper bound, and since $n\ge 5$, we also have
\begin{equation*}
    \frac{n (n-1) (n-2)}{32 \cdot \mathds{1}_{n=5}+(n^2+3n)\mathds{1}_{n\ge 6}} = \frac{60}{32}\mathds{1}_{n=5} + \frac{n (n^2-3n)}{(n-1)(n+2)}\mathds{1}_{n\ge 6} >
    \\
    (n-4). 
\end{equation*}
Therefore, $\boxi(G) \geq n-3$, which proves the corollary.
\end{proof}

\subsection{Plan of the paper}
In Section~\ref{sec prelims} we introduce several preliminary results. In Section~\ref{sec upper bound} we prove the upper bound in Theorem~\ref{main thm}. In Section~\ref{sec lower bounds}  we prove the lower bound in Theorem~\ref{main thm} and we discuss how to obtain general lower bounds through some already known result about poset dimension. In Section~\ref{sec 5} we prove Theorem \ref{thm complements line graphs}. We conclude the paper with a related discussion in Section~\ref{sec: conclusion}.

\section{Preliminaries}\label{sec prelims}

\subsection{Preliminaries on interval graphs}
Let $V$ be a ground set and $\mathcal{F}$ a family of subsets of $V$. The \emph{intersection graph} of $\mathcal{F}$ is the graph with vertex set $\mathcal{F}$ and edge set $\{S_1S_2:\hspace{0.2em} S_1, S_2\in \mathcal F, S_1\cap S_2\neq \varnothing\}$ 
A graph $G$ is an \emph{interval graph} if it can be represented as the intersection graph of a family of closed subintervals of the real line (see  Golumbic \cite{2004_Golumbic} or Gy\'arf\'as \cite{2003_Gyarfas} for a survey).

In view of the proof of the upper bound in Theorem \ref{main thm}, it will be useful to restate the geometric definition of boxicity in terms of interval graphs.

\begin{observation}
[\cite{1983_Cozzens}, Theorem 3] \label{box equi}
A graph $G = (V,E)$ has boxicity at most $k$ if and only if there are $k$ interval graphs $I_i = (V, E_i)$ for $i \in [k]$ such that $E(G) = \bigcap_{i \in [k]} E(I_i)$, or equivalently $E(G^c) = \bigcup_{i \in [k]} E(I^c_i)$.
\end{observation}

\subsection{Graph theoretic preliminaries}
For a graph $G = (V,E)$ and a set $S\subseteq V$, denote 
\begin{equation*}
    N_G(S) = \{u\in V\setminus S\hspace{0.2em}|\hspace{0.2em} \forall v\in S, uv\in E\} .
\end{equation*}
Moreover, let
\begin{equation*}
    c(k, G) = \max_{S\subseteq V: |S| = k} |N_G(S)|.
\end{equation*}

Note that for any $i,j\in \mathbb N$, $c(i,G) = j$ if and only if $G$ contains a copy of the complete bipartite graph $K_{i,j}$ as a subgraph, and contains no copy of $K_{i,j+1}$.

The following lemma appears as Theorem 2 in \cite{ACS}.
\begin{lemma}[\cite{ACS}, Theorem 2]\label{ACS lemma}
Let $G$ be a non-complete graph on $n$ vertices. Then,
\begin{equation*}
    \mathrm{box}(G)\ge \dfrac{|E(G^c)|}{\sum_{i=1}^{n-1} c(i, G^c)}.
\end{equation*}
\end{lemma}

\subsection{Other preliminaries}
We finish the preliminary section with two classical inequalities.

\begin{lemma}[Bernoulli's inequality, see e.g. \cite{Car}]\label{lem Bernoulli}
For every real number $a > -1$ and for every positive integer $n$ we have that $(1+a)^n\ge 1+an$.
\end{lemma}

\noindent
We finally state the Erd\H{o}s-Ko-Rado Theorem~\cite{EKR} - one of the most fundamental results in set theory.

\begin{theorem}[\cite{EKR}]\label{EKR thm}
Fix two positive integers $k,n$ with $2k\le n$. Let $\mathcal A$ be a family of $k$-subsets of $[n]$ such that any pair of sets have a non-empty intersection. Then, $|\mathcal A|\le \binom{n-1}{k-1}$, and if $2k+1\le n$, equality holds only for families of $k$-subsets of $[n]$, all containing a fixed element $i\le n$.
\end{theorem}

\section{Proof of the upper bound in Theorem~\ref{main thm}}\label{sec upper bound}
Fix two positive integers $k,n$ with $4\le 2k\le n-1$. We will construct a covering of the complement of the graph $Kn(k,n)$ with complements of $n-2$ interval graphs, and then conclude by Theorem \ref{box equi}. Below, we adopt the convention that $\binom{a}{b} = 0$ if $b < 0$.
For every $i\in [n-2]$, define the interval graph $I_i$ as follows.
\begin{itemize}
    \item To every set $S\subseteq [n], |S|=k$, if neither of $i, n-1, n$ is in $S$, assign the interval $\mathbb R$ to $S$.
    \item Assign each of the intervals 
    \begin{equation*}
        [2i, 2i+1]_{0\le i\le \binom{n-3}{k-3} - 1}
    \end{equation*}
    to a different set among the $\binom{n-3}{k-3}$ $k$-subsets of $[n]$ containing $\{i, n-1, n\}$.
    
    \item Assign each of the intervals 
    \begin{equation*}
    \left[2i+2\binom{n-3}{k-3}, 2i+2\binom{n-3}{k-3}+1\right]_{0\le i\le \binom{n-3}{k-2} - 1}
    \end{equation*}
    to a different set among the $\binom{n-3}{k-2}$ $k$-subsets of $[n]$, containing $i$ and $n$, but not $n-1$.
    \item Assign each of the intervals
    \begin{equation*}
    \left[2i+2\binom{n-3}{k-3}+2\binom{n-3}{k-2}, 2i+2\binom{n-3}{k-3}+2\binom{n-3}{k-2}+1\right]_{0\le i\le \binom{n-3}{k-1} - 1} 
    \end{equation*}
    to a different set among the $\binom{n-3}{k-1}$ $k$-subsets of $[n]$, containing $i$, but neither $n$ nor $n-1$.
    
    \item Assign each of the intervals
    \begin{equation*}
    \left[2i+2\binom{n-3}{k-3}+2\binom{n-3}{k-2}+2\binom{n-3}{k-1}, 2i+2\binom{n-3}{k-3}+2\binom{n-3}{k-2}+2\binom{n-3}{k-1}+1\right]_{0\le i\le \binom{n-3}{k-2} - 1}   
    \end{equation*}
    to a different set among the $\binom{n-3}{k-2}$ $k$-subsets of $[n]$, containing $i$ and $n-1$, but not $n$.
    
    \item Assign the interval
    \begin{equation*}
    \left[2\binom{n-3}{k-3}+2\binom{n-3}{k-2}, 2\binom{n-3}{k-3}+2\binom{n-3}{k-2}+2\binom{n-3}{k-1}-1\right]
    \end{equation*}
    to all of the $\binom{n-3}{k-2}$ $k$-subsets of $[n]$, containing $n-1$ and $n$, but not $i$.
    
    \item Assign the interval
    \begin{equation*}
    \left[2\binom{n-3}{k-3}, 2\binom{n-3}{k-3}+2\binom{n-3}{k-2}+2\binom{n-3}{k-1}-1\right]
    \end{equation*}
    to all of the $\binom{n-3}{k-1}$ $k$-subsets of $[n]$, containing $n-1$, but neither $n$ nor $i$.
    
    \item Assign each of the intervals
    \begin{equation*}
    \left[2\binom{n-3}{k-3}+2\binom{n-3}{k-2}, 2\binom{n-3}{k-3}+4\binom{n-3}{k-2}+2\binom{n-3}{k-1}-1\right]  
    \end{equation*}
    to a different set among the $\binom{n-3}{k-1}$ $k$-subsets of $[n]$, containing $n$, but neither $n-1$ nor $i$.
\end{itemize}

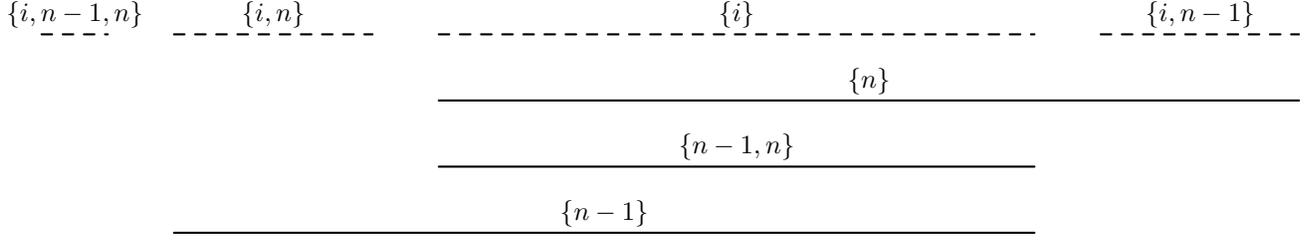
\begin{figure}
\centering
\begin{tikzpicture}[scale=0.88,line cap=round,line join=round,x=1cm,y=1cm]
\clip(-10,-1.5) rectangle (10,2.5);
\draw [line width=0.8pt,dash pattern=on 4pt off 4pt] (-9,2)-- (-8,2);
\draw [fill=black] (-8.5,2.3) node {$\{i,n-1,n\}$};

\draw [line width=0.8pt,dash pattern=on 4pt off 4pt] (-7,2)-- (-4,2);
\draw [fill=black] (-5.5,2.3) node {$\{i,n\}$};

\draw [line width=0.8pt,dash pattern=on 4pt off 4pt] (-3,2)-- (6,2);
\draw [fill=black] (1.5,2.3) node {$\{i\}$};

\draw [line width=0.8pt,dash pattern=on 4pt off 4pt] (7,2)-- (10,2);
\draw [fill=black] (8.5,2.3) node {$\{i,n-1\}$};

\draw [line width=0.8pt] (-3,1)-- (10,1);
\draw [fill=black] (3.5,1.3) node {$\{n\}$};

\draw [line width=0.8pt] (-3,0)-- (6,0);
\draw [fill=black] (1.5,0.3) node {$\{n-1,n\}$};

\draw [line width=0.8pt] (-7,-1)-- (6,-1);
\draw [fill=black] (-0.5,-0.7) node {$\{n-1\}$};

\end{tikzpicture}
\caption{For $i\in [n-2]$, the figure represents the positions of all finite intervals in $I_i$. Dotted lines correspond to a number of consecutive disjoint intervals, solid lines correspond to single intervals. On every line (solid or dotted) is denoted the exact subset of $\{i, n-1, n\}$, which is included in the sets, corresponding to the particular interval or group of disjoint intervals.}
\label{fig 1}
\end{figure}
\noindent
Figure~\ref{fig 1} shows a representation of the described intervals.

For every $i\in [n]$, denote by $K_i$ the complete graph on all vertices in $Kn(k,n)$, corresponding to sets, containing $i$. One may readily check that:
\begin{itemize}
    \item $\bigcup_{j \in [n]} E(K_j) = E(Kn(k,n)^c)$,
    \item for every $i\in [n-2]$, $E(K_i^c) \subset E(I^c_i) \subseteq E(Kn(k,n)^c)$, and
    \item the cliques $K_{n-1}$ and $K_n$ are both contained in $\cup_{i\in [n-2]} I^c_i$.
\end{itemize}
This shows that $\cup_{i\in [n-2]} E(I^c_i) = E(Kn(k, n)^c)$, concluding the proof.
\qed

\section{Proof of the lower bound in Theorem~\ref{main thm}}\label{sec lower bounds}

For any $k\ge 2$ and $n\ge 2k^3-2k^2+1$, fix $G = Kn(k,n)$ and $c(\cdot) = c(\cdot, G^c)$, where $c(\cdot, \cdot)$ was defined just before Lemma~\ref{ACS lemma}.

\begin{observation}\label{ob trivial}
For any $i,j\in \mathbb N$, if $c(i) = j$, then $c(j+1)< i\le c(j)$.
\end{observation}
\begin{proof}
The fact that $c(i) = j$ means that there is a copy of the complete bipartite graph $K_{i,j}$, included in $G^c$, but no copy of $K_{i, j+1}$ could be realised as a subgraph of $G^c$, which is equivalent to our claim.
\end{proof}

\begin{corollary}\label{cor trivial}
For any $i,j\in \mathbb N$, if $c(i) = j$ and $c(i+1) < j$, then $c(j) = i$.
\end{corollary}
\begin{proof}
Fix $s = c(i+1)$. By Observation~\ref{ob trivial} for $(i,j)$ we have that $i\le c(j)$, and since $s < j$, by the same result for $(i+1,j)$ we deduce that $c(j)\le c(s+1) < i+1 \le c(s)$, which proves the corollary.
\end{proof}

\begin{figure}
\centering
\begin{tikzpicture}[scale=0.6, line cap=round,line join=round,x=1cm,y=1cm]
\clip(-16,-3.58) rectangle (12.8,7.5);
\draw [->,line width=0.5pt] (-6,-3) -- (4,-3);
\draw [->,line width=0.5pt] (-6,-3) -- (-6,7);
\draw [line width=0.5pt] (-6,6)-- (-5,6);
\draw [line width=0.5pt] (-5,6)-- (-5,-3);
\draw [line width=0.5pt] (3,-3)-- (3,-2);
\draw [line width=0.5pt] (3,-2)-- (-6,-2);
\draw [line width=0.5pt] (-6,4)-- (-4,4);
\draw [line width=0.5pt] (-4,4)-- (-4,-3);
\draw [line width=0.5pt] (1,-3)-- (1,-1);
\draw [line width=0.5pt] (1,-1)-- (-6,-1);
\draw [line width=0.5pt] (-6,3)-- (-3,3);
\draw [line width=0.5pt] (-3,3)-- (-3,-3);
\draw [line width=0.5pt] (0,-3)-- (0,0);
\draw [line width=0.5pt] (0,0)-- (-6,0);
\draw [line width=0.5pt] (-6,2)-- (-1,2);
\draw [line width=0.5pt] (-1,2)-- (-1,-3);
\draw [line width=0.5pt] (-1,1)-- (-6,1);
\draw [line width=0.5pt] (-2,2)-- (-2,-3);
\draw [line width=0.5pt] (2,-2)-- (2,-3);
\draw [line width=0.5pt] (-6,5)-- (-5,5);
\draw [line width=0.5pt] (-6,-3)-- (4,7);
\begin{scriptsize}
\draw [fill=black] (4.3,-3.2) node {\large{$i$}};
\draw [fill=black] (-7,7) node {\large{$c(i)$}};
\draw [fill=black] (-7,-3.2) node {\large{$(0,0)$}};
\end{scriptsize}
\end{tikzpicture}
\caption{A picture of the Young diagram, coming from the sequence $(9,7,6,5,5,3,2,1,1)$. Here, $t = 5$.}
\label{fig 2}
\end{figure}
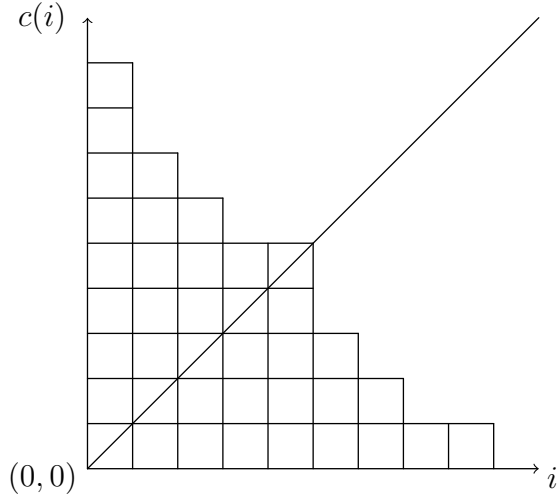

A more visual interpretation of the last corollary is the following. Consider a Young diagram with columns of altitude $(c(i))_{1\le i\le \binom{n}{k}-1}$, that is, for every $i\in \left[\binom{n}{k}-1\right]$, the column over $[i-1,i]$ has height $c(i)$ (see Figure~\ref{fig 2}). Then, Observation~\ref{ob trivial} and Corollary~\ref{cor trivial} imply that this Young diagram is symmetric with respect to the line $y=x$. It follows that, when computing $\sum_{i=1}^{\binom{n}{k}-1} c(i)$, it is sufficient to compute the area of the diagram above the line $y=x$ and to multiply by two. Let $t = \max\{i: c(i)\ge i\}$. The expression of this area is given by
\begin{equation}\label{eq area}
    \sum_{i=1}^{t} \big (c(i) - (i-1/2) \big) = \sum_{i=1}^t c(i) - \dfrac{t^2}{2}.
\end{equation}

Call a bipartite graph with given parts $(V_1, V_2)$ \emph{balanced} if $|V_1| = |V_2|$. 

\begin{lemma}\label{lem bipartite}
The largest balanced complete bipartite graph, contained in $G^c$ as a subgraph, contains $2\left\lfloor \frac{1}{2}\binom{n-1}{k-1}\right\rfloor$ vertices, or equivalently $t = \left\lfloor \frac{1}{2}\binom{n-1}{k-1}\right\rfloor$.
\end{lemma}
\begin{proof}
Denote by $s$ the number of vertices in every part of a largest balanced complete bipartite subgraph of $G^c$. First of all, there is a clique in $G^c$ that contains $\binom{n-1}{k-1}$ vertices, so $s\ge \left\lfloor \frac{1}{2}\binom{n-1}{k-1}\right\rfloor$. 

To prove the upper bound, we argue by contradiction. Suppose that $2s \ge \binom{n-1}{k-1}+1$. Then, by Theorem~\ref{EKR thm} there are two vertices, corresponding to disjoint $k$-subsets $A$ and $B$ of $[n]$, and clearly these vertices must be contained in the same part. Then, each of the sets corresponding to a vertex in the other part must contain one element from both $A$ and $B$. Thus, there are at most $k^2 \binom{n-2}{k-2}$ vertices in every part, and since $n\ge 2k^3 - 2k^2 + 1$, one may deduce that
\begin{equation*}
    2s\le 2k^2\binom{n-2}{k-2}\le \binom{n-1}{k-1}.
\end{equation*}
This contradiction concludes the proof of the lemma.
\end{proof}

\begin{observation}\label{cor Pascal}
For every triplet of positive integers $(a,b,c)$ such that $a>c$,
\begin{equation*}
    \binom{a}{b} - \binom{a-c}{b} \leq c \binom{a-1}{b-1}.
\end{equation*}
\begin{proof}
Using Pascal’s identity $c$ times we deduce that 
\begin{equation*}
    \binom{a}{b} - \binom{a-c}{b} = \sum_{i=1}^{c} \binom{a-i+1}{b} - \binom{a-i}{b} = \sum_{i=1}^{c} \binom{a-i}{b-1}.
\end{equation*}
Then, the result follows by the trivial inequality $\binom{a-i}{b-1}\le \binom{a-1}{b-1}$ for every $i\ge 1$.
\end{proof}
\end{observation}

\begin{corollary}\label{range 0}
For every $k\ge 2$ and $n\ge 2k^3-2k^2+1$,
\begin{equation*}
    \sum_{i=1}^{\binom{n-2}{k-2}-\binom{n-2-k}{k-2}} c(i)\le  k^2\binom{n-3}{k-3}\binom{n-1}{k-1}.
\end{equation*}
\end{corollary}
\begin{proof}
For every $i$ we have that $c(i)\le c(1) = \binom{n}{k}-\binom{n-k}{k}-1$, so
\begin{equation*}
    \sum_{i=1}^{\binom{n-2}{k-2}-\binom{n-2-k}{k-2}} c(i)\le \left(\binom{n-2}{k-2}-\binom{n-2-k}{k-2}\right)\left(\binom{n}{k}-\binom{n-k}{k}-1\right)\le k^2\binom{n-3}{k-3}\binom{n-1}{k-1}.
\end{equation*}
where the last inequality is achieved by two consecutive applications of Observation~\ref{cor Pascal}.
\end{proof}

Fix a family $\mathcal F_i$ of $i$ $k$-subsets of $[n]$. 

\begin{lemma}\label{1-intersection}
Suppose that the intersection of all sets in the family $\mathcal F_i$ contains exactly one element of $[n]$. Then, there are at most $\binom{n-1}{k-1} - i + (k-1)^2\binom{n-3}{k-2}$ $k$-subsets of $[n]$ that are not contained in $\mathcal F_i$ and intersect each of the sets in $\mathcal F_i$.
\end{lemma}
\begin{proof}
Without loss of generality let the common element of all sets in $\mathcal F_i$ be 1. Let $A = \{a_1 = 1, a_2, \dots, a_k\}$ be a member of $\mathcal F_i$, and for every $j\in [k]\setminus 1$, let $B_j$ be a set in $\mathcal F_i$, not containing $a_j$. Then, any $k$-subset of $[n]$ which intersects all members of $\mathcal F_i$ either contains 1 or it contains an element $a_j$ among $\{a_2, \dots, a_k\}$, and at least one of the $k-1$ elements of $B_j$ different from 1. Thus, there are
\begin{equation*} 
    \binom{n-1}{k-1} - i
\end{equation*} 
$k$-sets outside $\mathcal F_i$, containing 1, and at most
\begin{equation*}
    (k-1)^2\binom{n-3}{k-2}
\end{equation*}
$k$-subsets of $[n]$, which intersect every element of $\mathcal F_i$, but do not contain 1. This proves the lemma.
\end{proof}

\begin{lemma}\label{0-intersection}
Suppose that the intersection of all sets in the family $\mathcal F_i$ is empty. Then, there are at most $k^2\binom{n-2}{k-2}$ that are not contained in $\mathcal F_i$ and intersect each of the sets in $\mathcal F_i$.
\end{lemma}
\begin{proof}
Let $A = \{a_1, \dots, a_k\}$ be an arbitrary set in $\mathcal F_i$, and for every $j\in [k]$, let $B_j$ be an arbitrary set in $\mathcal F_i$, not containing $a_j$. Counting all $k$-subsets, intersecting all members of $\mathcal F_i$, according to the first element in $A$ which they contain, say $a_j$, and then the first element of $B_j$ which they contain, we get at most $k^2\binom{n-2}{k-2}$
$k$-subsets of $[n]$ outside $\mathcal F_i$, which intersect each of the $k$-sets in $\mathcal F_i$.
\end{proof}

\begin{lemma}\label{range 1}
For every $k\ge 2$, $n\ge 2k^3-2k^2+1$ and $i\in \left[\binom{n-2}{k-2}-\binom{n-2-k}{k-2}+1, \binom{n-2}{k-2}\right]$, $c(i) = 2\binom{n-1}{k-1} - \binom{n-2}{k-2}-i$.
\end{lemma}
\begin{proof}
We prove that, whatever the choice of $\mathcal F_i$, the vertices in $G^c$, corresponding to the sets in $\mathcal F_i$, have at most $2\binom{n-1}{k-1} - \binom{n-2}{k-2}-i$ common neighbours. Since $i > \binom{n-3}{k-3}$, the sets in $\mathcal F_i$ could have at most two common elements. If the given $i$ sets contain two common elements, then let these elements be 1 and 2. We show that there is no $k$-subset $\{a_1,a_2,\dots,a_k\}\subseteq [n]\setminus \{1,2\}$ intersecting all the members of $\mathcal{F}_i$.
Indeed, the number of all $k$-subsets of $[n]$, containing 1, 2 and at least one element among $\{a_1,a_2,\dots,a_k\}$, is $\binom{n-2}{k-2}-\binom{n-2-k}{k-2} < i = |\mathcal{F}_i|$.
Thus, any $k$-set, intersecting all members of $\mathcal F_i$, contains either 1 or 2, or both, and there are exactly $2\binom{n-1}{k-1} - \binom{n-2}{k-2}-i$ such subsets of $[n]$ outside $\mathcal F_i$.\par

By Lemma~\ref{1-intersection} and Lemma~\ref{0-intersection} it remains to verify that for every $i\in \left[\binom{n-2}{k-2}-\binom{n-2-k}{k-2}+1, \binom{n-2}{k-2}\right]$
\begin{equation*}
    \max\left(k^2\binom{n-2}{k-2}, \binom{n-1}{k-1} - i + (k-1)^2\binom{n-3}{k-2}\right)\le 2\binom{n-1}{k-1} - \binom{n-2}{k-2}-i,
\end{equation*}
which is equivalent to
\begin{equation*}
    k^2\binom{n-2}{k-2}\le 2\binom{n-1}{k-1} - 2\binom{n-2}{k-2} = 2\binom{n-2}{k-1}
    \text{ and }
    (k-1)^2\binom{n-3}{k-2}\le \binom{n-1}{k-1} - \binom{n-2}{k-2} = \binom{n-2}{k-1}.
\end{equation*}

Since $n\ge k^3$, we have
\begin{align*}
    2\binom{n-2}{k-1} = 2\dfrac{n-k}{k-1}\binom{n-2}{k-2}\ge k^2\binom{n-2}{k-2},
\end{align*}
and the second inequality holds since
\begin{equation*}
    \binom{n-2}{k-1} = \dfrac{n-2}{k-1}\binom{n-3}{k-2}\ge (k-1)^2\binom{n-3}{k-2}.
\end{equation*}

The lemma is proved.
\end{proof}

\begin{lemma}\label{range 2}
For every $k\ge 2$, $n\ge 2k^3-2k^2+1$ and $i\in \left[\binom{n-2}{k-2}+1, \binom{n-1}{k-1}-\binom{n-1-k}{k-1}\right]$, $c(i) = \binom{n-1}{k-1} - i + (k-1)^2\binom{n-3}{k-2}$.
\end{lemma}
\begin{proof}
Once again, we work with $\mathcal F_i$. Since $i > \binom{n-2}{k-2}$, the sets in $\mathcal F_i$ could have at most one common element. Thus, by Lemma~\ref{1-intersection} and Lemma~\ref{0-intersection} it remains to prove that for every $i\le \binom{n-1}{k-1}-\binom{n-1-k}{k-1}$ we have
\begin{equation*}
\binom{n-1}{k-1} - i + (k-1)^2\binom{n-3}{k-2}\ge k^2\binom{n-2}{k-2}.
\end{equation*}

Note that, on the one hand,
\begin{equation*}
\binom{n-1-k}{k-1} = \left(\prod_{j=0}^{k-2} \dfrac{n-1-j-k}{n-1-j}\right)\binom{n-1}{k-1}\ge \left(1 - \dfrac{k}{n-k+1}\right)^{k-1}\dfrac{n-1}{k-1}\binom{n-2}{k-2},
\end{equation*}
and on the other hand, by Bernoulli's inequality (Lemma~\ref{lem Bernoulli}) and the assumption $n\ge 2k^3-2k^2+1$ we have
\begin{equation*}
\left(1 - \dfrac{k}{n-k+1}\right)^{k-1}\dfrac{n-1}{k-1}\ge \left(1-\dfrac{k(k-1)}{n-k+1}\right) 2k^2\ge \left(1-\dfrac{k(k-1)}{2(k-1)^2(k+1)}\right) 2k^2\ge k^2.
\end{equation*}
\end{proof}

\begin{lemma}\label{range 3}
For every $k\ge 2$, $n\ge 2k^3-2k^2+1$ and $i\in \left[\binom{n-1}{k-1}-\binom{n-1-k}{k-1}+1, \left\lfloor\frac{1}{2}\binom{n-1}{k-1}\right\rfloor\right]$, $c(i) = \binom{n-1}{k-1}-i$.
\end{lemma}
\begin{proof}
Once again, we work with $\mathcal F_i$. Since $i > \binom{n-2}{k-2}$, the sets in $\mathcal F_i$ could have at most one common element. We consider two cases.\par

If the sets in $\mathcal F_i$ all contain one common element, let this element be 1 without loss of generality. Then, since $i\ge \binom{n-1}{k-1}-\binom{n-1-k}{k-1}+1$, any $k$-subset of $[n]$, intersecting all members in the family, contains 1. Indeed, the number of all $k$-subsets of $[n]$, containing 1 and containing some element among $\{a_1,a_2,\dots,a_k\}$ for any $\{a_1,a_2,\dots,a_k\}\subseteq [n]\setminus 1$, is $\binom{n-1}{k-1}-\binom{n-1-k}{k-1} < i$. Since the number of $k$-sets, containing 1 and not in $\mathcal F_i$, is $\binom{n-1}{k-1}-i$, we conclude that $c(i)\ge \binom{n-1}{k-1}-i$.\par

If the sets in $\mathcal F_i$ do not have a common element, by Lemma~\ref{0-intersection} there are at most $k^2\binom{n-2}{k-2}$ elements outside $\mathcal F_i$, which intersect each of the $k$-sets in $\mathcal F_i$. It remains to observe that, for $n\ge 2k^3-2k^2+1$, $k^2\binom{n-2}{k-2}\le \frac{1}{2}\binom{n-1}{k-1}$ and so $k^2\binom{n-2}{k-2}\le \binom{n-1}{k-1} - \left\lfloor\frac{1}{2}\binom{n-1}{k-1}\right\rfloor\le \binom{n-1}{k-1} - i$, which proves the lemma.
\end{proof}

Using \eqref{eq area} we deduce that
\begin{equation*}
    \dfrac{1}{2}\sum_{i=1}^{\binom{n}{k}-1} c(i) \le\hspace{0.3em} \sum_{i=1}^{\lfloor \frac{1}{2}\binom{n-1}{k-1}\rfloor} c(i) - \dfrac{1}{2}\left\lfloor \frac{1}{2}\binom{n-1}{k-1}\right\rfloor^2.
\end{equation*}

We separate the above sum into four sums over the intervals $\left[1, \binom{n-2}{k-2}-\binom{n-2-k}{k-2}\right]$, $\left[\binom{n-2}{k-2}-\binom{n-2-k}{k-2}+1, \binom{n-2}{k-2}\right]$, $\left[\binom{n-2}{k-2}+1, \binom{n-1}{k-1}-\binom{n-1-k}{k-1}\right]$ and $\left[\binom{n-1}{k-1}-\binom{n-1-k}{k-1}+1, \lfloor \frac{1}{2}\binom{n-1}{k-1}\rfloor\right]$. Then, by Lemma~\ref{range 0} we get that
\begin{equation*}
    \sum_{i=1}^{\binom{n-2}{k-2}-\binom{n-2-k}{k-2}} c(i)\le k^2 \binom{n-3}{k-3}\binom{n-1}{k-1} = \dfrac{k^2(k-1)(k-2)}{(n-1)(n-2)} \binom{n-1}{k-1}^2\le \dfrac{k^2(k-1)^2}{(n-1)^2}\binom{n-1}{k-1}^2,
\end{equation*}
where the last inequality holds since $\dfrac{k-2}{n-2}\le \dfrac{k-1}{n-1}$.

By Lemma~\ref{range 1} we get that
\begin{align*}
    \sum_{i=\binom{n-2}{k-2}-\binom{n-2-k}{k-2}+1}^{\binom{n-2}{k-2}} c(i)
    &\le\hspace{0.3em} \sum_{i=\binom{n-2}{k-2}-\binom{n-2-k}{k-2}+1}^{\binom{n-2}{k-2}} \left(2\binom{n-1}{k-1} - \binom{n-2}{k-2}-i\right)\\
    &\le\hspace{0.3em} \binom{n-2-k}{k-2} \left(2\binom{n-1}{k-1} - \binom{n-2}{k-2}\right)-\sum_{i=\binom{n-2}{k-2}-\binom{n-2-k}{k-2}+1}^{\binom{n-2}{k-2}} i\\
    &\le\hspace{0.3em} \binom{n-2}{k-2} \left(2-\dfrac{k-1}{n-1}\right)\binom{n-1}{k-1}-\sum_{i=\binom{n-2}{k-2}-\binom{n-2-k}{k-2}+1}^{\binom{n-2}{k-2}} i\\
    &\le\hspace{0.3em} \dfrac{(k-1)(2n-k-1)}{(n-1)^2}\binom{n-1}{k-1}^2-\sum_{i=\binom{n-2}{k-2}-\binom{n-2-k}{k-2}+1}^{\binom{n-2}{k-2}} i\\
    &\le\hspace{0.3em} \dfrac{2(k-1)}{n-1}\binom{n-1}{k-1}^2-\sum_{i=\binom{n-2}{k-2}-\binom{n-2-k}{k-2}+1}^{\binom{n-2}{k-2}} i.
\end{align*}

By Lemma~\ref{range 2} we get that
\begin{align*}
    \sum_{i=\binom{n-2}{k-2}+1}^{\binom{n-1}{k-1}-\binom{n-1-k}{k-1}} c(i)
    &\le\hspace{0.3em} \sum_{i=\binom{n-2}{k-2}+1}^{\binom{n-1}{k-1}-\binom{n-1-k}{k-1}} \left(\binom{n-1}{k-1} - i + (k-1)^2\binom{n-3}{k-2}\right)\\
    &\le\hspace{0.3em} \left(\binom{n-1}{k-1} - \binom{n-1-k}{k-1} - \binom{n-2}{k-2}\right) \left(\binom{n-1}{k-1} + (k-1)^2\binom{n-3}{k-2}\right) - \sum_{i=\binom{n-2}{k-2}+1}^{\binom{n-1}{k-1}-\binom{n-1-k}{k-1}} i\\
    &\le\hspace{0.3em} (k-1)\binom{n-2}{k-2}\left(1+\dfrac{(k-1)^3}{n-1}\right)\binom{n-1}{k-1} - \sum_{i=\binom{n-2}{k-2}+1}^{\binom{n-1}{k-1}-\binom{n-1-k}{k-1}} i\\
    &=\hspace{0.3em} \dfrac{(k-1)^2}{n-1}\left(1+\dfrac{(k-1)^3}{n-1}\right)\binom{n-1}{k-1}^2 - \sum_{i=\binom{n-2}{k-2}+1}^{\binom{n-1}{k-1}-\binom{n-1-k}{k-1}} i.
\end{align*}

Also, by Lemma~\ref{range 3} we have
\begin{align*}
\sum_{i=\binom{n-1}{k-1}-\binom{n-1-k}{k-1}+1}^{\lfloor \frac{1}{2}\binom{n-1}{k-1}\rfloor} c(i)
&\le\hspace{0.3em} \sum_{i=\binom{n-1}{k-1}-\binom{n-1-k}{k-1}+1}^{\lfloor \frac{1}{2}\binom{n-1}{k-1}\rfloor}\left(\binom{n-1}{k-1} - i\right)\\
&\le\hspace{0.3em} \left(\binom{n-1-k}{k-1} - \left\lceil \frac{1}{2}\binom{n-1}{k-1}\right\rceil \right)\binom{n-1}{k-1} - \sum_{i=\binom{n-1}{k-1}-\binom{n-1-k}{k-1}+1}^{\lfloor \frac{1}{2}\binom{n-1}{k-1}\rfloor} i\\
&\le\hspace{0.3em} \dfrac{1}{2}\binom{n-1}{k-1}^2 - \sum_{i=\binom{n-1}{k-1}-\binom{n-1-k}{k-1}+1}^{\lfloor \frac{1}{2}\binom{n-1}{k-1}\rfloor} i.
\end{align*}

Focusing on the sums which are not yet developed, we have that 
\begin{align*}
&\left(\sum_{i=\binom{n-2}{k-2}-\binom{n-2-k}{k-2}+1}^{\lfloor \frac{1}{2}\binom{n-1}{k-1}\rfloor} i\right) + \frac{1}{2}\left\lfloor \frac{1}{2}\binom{n-1}{k-1}\right\rfloor^2\\
\ge\hspace{0.3em} 
&\left(\sum_{i=1}^{\lfloor \frac{1}{2}\binom{n-1}{k-1}\rfloor} i\right) - \left(\sum_{i=1}^{\binom{n-2}{k-2}-\binom{n-2-k}{k-2}} i\right) + \frac{1}{2}\left\lfloor \frac{1}{2}\binom{n-1}{k-1}\right\rfloor^2
\\
\ge\hspace{0.3em} 
&\dfrac{1}{2}\left\lfloor \frac{1}{2}\binom{n-1}{k-1}\right\rfloor^2 +\dfrac{1}{2}\left\lfloor \frac{1}{2}\binom{n-1}{k-1}\right\rfloor - \dfrac{1}{2}\left(\binom{n-2}{k-2}-\binom{n-2-k}{k-2} + 1\right)^2+\dfrac{1}{2}\left\lfloor \frac{1}{2}\binom{n-1}{k-1}\right\rfloor^2
\\
\ge\hspace{0.3em} 
&\left( \frac{1}{2}\binom{n-1}{k-1} - \frac{1}{2}\right)^2 +\dfrac{1}{2}\left( \frac{1}{2}\binom{n-1}{k-1} - \frac{1}{2}\right) - \dfrac{1}{2}\left(\binom{n-2}{k-2}-\binom{n-2-k}{k-2} + 1\right)^2
\\
\ge\hspace{0.3em} 
&\dfrac{1}{4}\binom{n-1}{k-1}^2 - \frac{1}{4}\binom{n-1}{k-1} - \dfrac{1}{2}\left(\binom{n-2}{k-2}-\binom{n-2-k}{k-2} + 1\right)^2
\\
\ge\hspace{0.3em} 
&\dfrac{1}{4}\binom{n-1}{k-1}^2 - \frac{1}{4}\binom{n-1}{k-1} - \dfrac{1}{2}\binom{n-2}{k-2}^2
\\
\ge\hspace{0.3em} 
&\left(\dfrac{1}{4} - \dfrac{k-1}{4(n-1)} - \dfrac{(k-1)^2}{2(n-1)^2}\right)\binom{n-1}{k-1}^2.
\end{align*}

We conclude that $\sum_{i=1}^{\binom{n}{k} - 1} c(i)$ is bounded from below by
\begin{align}\label{eq bound 1}
& 2 \left(\dfrac{k^2(k-1)^2}{(n-1)^2} + \dfrac{2(k-1)}{n-1} + \dfrac{(k-1)^2}{n-1}\left(1+\dfrac{(k-1)^3}{n-1}\right) + \dfrac{1}{2} - \left(\dfrac{1}{4}-\dfrac{(k-1)}{4(n-1)}-\dfrac{(k-1)^2}{2(n-1)^2}\right)\right) \binom{n-1}{k-1}^2\nonumber
\\
=\hspace{0.3em}
& 2\left(\dfrac{1}{4} + \dfrac{(k-1)^2 + 9(k-1)/4}{n-1} + \dfrac{(k-1)^5 + k^2 (k-1)^2 + (k-1)^2/2}{(n-1)^2}\right) \binom{n-1}{k-1}^2\nonumber
\\
\le\hspace{0.3em}
& 2\left(\dfrac{1}{4} + \dfrac{(k-1)^2 + 9(k-1)/4}{n-1} + \dfrac{((k-1)^5 + (k-1)^2/2) + k^2 (k-1)^2}{2k^2(k-1)(n-1)}\right) \binom{n-1}{k-1}^2\nonumber
\\
\le\hspace{0.3em}
& \left(\dfrac{1}{2} + \dfrac{2(k-1)^2 + 9(k-1)/2}{n-1} + \dfrac{(k^2-4k+6) + k-1}{n-1}\right) \binom{n-1}{k-1}^2\nonumber
\\
\le\hspace{0.3em}
& \left(\dfrac{1}{2} + \dfrac{3k^2 - 5k/2 + 5/2}{n-1}\right) \binom{n-1}{k-1}^2.
\end{align}

We conclude by \eqref{eq bound 1} and Lemma~\ref{ACS lemma} that $\boxi(G)$ is at least
\begin{equation*}
\dfrac{|E(G^c)|}{\left(\dfrac{1}{2} + \dfrac{3k^2-5k/2+5/2}{n-1}\right) \binom{n-1}{k-1}^2}\ge \dfrac{\dfrac{1}{2}\binom{n}{k}\left(\binom{n}{k} - \binom{n-k}{k}\right)}{\left(\dfrac{1}{2} + \dfrac{3k^2-5k/2+5/2}{n-1}\right) \binom{n-1}{k-1}^2} = \dfrac{n}{\left(1+\dfrac{6k^2-5k+5}{n-1}\right)}\dfrac{\binom{n}{k} - \binom{n-k}{k}}{ k\binom{n-1}{k-1}}.
\end{equation*}

It remains to note that 
$$\dfrac{1}{1+\dfrac{6k^2-5k+5}{n-1}}\ge 1-\dfrac{6k^2-5k+5}{n-1}$$ 
and by Pascal's identity applied $k$ times
\begin{equation}\label{final eq}
    \binom{n}{k} - \binom{n-k}{k} = \sum_{i=1}^{k} \binom{n-i}{k-1} = \sum_{i=1}^{k} \binom{n-1}{k-1}\dfrac{\prod_{j=1}^{i-1} (n-k+1-j)}{\prod_{j=1}^{i-1} (n-j)}\ge \sum_{i=1}^{k} \binom{n-1}{k-1} \left(1 - \dfrac{k-1}{n-k+1}\right)^{i-1}.
\end{equation}
Since $n\ge 2k-1$, we have that $\frac{k-1}{n-k+1} < 1$, so by Lemma~\ref{lem Bernoulli} $\left(1 - \dfrac{k-1}{n-k+1}\right)^{i-1}\ge 1 - \dfrac{(k-1)(i-1)}{n-k+1}$, so we conclude that \eqref{final eq} is bounded from below by
\begin{equation*}
    \sum_{i=1}^{k} \binom{n-1}{k-1} \left(1 - \dfrac{(k-1)(i-1)}{n-k+1}\right)\ge \binom{n-1}{k-1} \left(k - \dfrac{(k-1)^2k}{2(n-k+1)}\right).
\end{equation*}

Thus, we get that $\boxi(G)$ is at least 
\begin{align*}
    \dfrac{n}{\left(1+\dfrac{6k^2-5k+5}{n-1}\right)}\dfrac{\binom{n}{k} - \binom{n-k}{k}}{ k\binom{n-1}{k-1}}
    &\ge\hspace{0.3em} n \left(1-\dfrac{6k^2-5k+5}{n-1}\right) \left(1 - \dfrac{(k-1)^2k}{2k(n-k+1)}\right)\\
    &\ge\hspace{0.3em} n - n\dfrac{6k^2-5k+5}{n-1} - \dfrac{n}{2}\dfrac{(k-1)^2}{n-k+1}\\
    &\ge\hspace{0.3em} n - n\dfrac{6k^2-5k+8}{n} - \dfrac{n}{2}\dfrac{(k-1)^2+k-1}{n}\\
    &\ge\hspace{0.3em} n - 6k^2+5k-8 - \dfrac{k(k-1)}{2}\\
    &=\hspace{0.3em} n - \dfrac{13k^2-11k+16}{2},
\end{align*}
where the third line comes from the inequalities
$$\dfrac{6k^2-5k+5}{n-1}\le 3 \text{ and } \dfrac{(k-1)^2}{(n-k+1)}\le 1,$$
which are ensured by the assumption that $n\ge 2k^3 - 2k^2+1$. The proof of Theorem~\ref{main thm} is finished.
\qed

\begin{remark}
In general, for any positive integers $k,n$ with $n\ge 2k+1$, lower bounds on the boxicity of the Kneser graph $Kn(k,n)$ can be easily derived thanks to the remarkable connection between graph boxicity and poset dimension, shown in~\cite{ABC}. To explain the approach, we define the \emph{extended double cover} of a graph $G$ to be the graph $G_c$ with vertex set $V_1\cup V_2$, where $V_1 = \{u_1:\hspace{0.2em} u\in V(G)\}$ and $V_2 = \{u_2:\hspace{0.2em} u\in V(G)\}$ are two disjoint copies of $V(G)$, and edge set $\{u_1v_2:\hspace{0.2em} u_1\in V_1, v_2\in V_2, u\equiv v \text{ or } uv\in E(G)\}$.

Fix $G = Kn(k,n)$, and let $S_u$ be the $k$-set which correspond to $u$ in the construction of $Kn(k,n)$. For every $u\in V(G)$, associate $S_u$ to $u_1\in V_1$  and $[n] \setminus S_u$ to $u_2 \in V_2$. Then, let $\Tilde{G}_c$ be the subgraph of $G_c$, induced by the vertices in $V_1\cup V_2$, whose corresponding sets contain the element 1. Clearly, $\boxi(\Tilde{G}_c)\le \boxi(G_c)$, and by Lemma~2 in \cite{ABC} we have $\boxi(G_c)-2\le \boxi(G)$.

Since $\Tilde{G}_c$ is the comparability graph of the poset $(k-1, n-k-1; n-1)$ (with elements the subsets of $[n-1]$ of size $k-1$ or $n-k-1$, partially ordered by the inclusion relation), combining Theorem~1 in~\cite{ACS} and Propositions~2.1~and~2.2 in~\cite{Fur}, one may deduce that
\begin{equation*}
    \boxi(G)\ge \frac{n-2k-1}{2},\ \text{if} \ 2k+1\le n\le 3k+1,
\end{equation*}
and 
\begin{equation*}
    \boxi(G)\ge \frac{n-k-4}{2},\ \text{if} \ n\ge 3k+2.
\end{equation*}
For values of $n$ ``close'' to $2k+1$, better lower bounds may be deduced by Theorem~4.5 in~\cite{Kie1}.
In particular, using crucially the fact that $(k-1, n-k-1; n-1)$ contains an isomorphic copy of $(1, n-2k+1; n-k+1)$ as an induced subposet, there is a positive constant $c > 0$ such that 
\begin{equation*}
    c 2^{n-2k+1} \log\log n - 2\le \boxi(G)\ \text{if} \ n-2k+1\le \log\log n - \log\log\log n.
\end{equation*}
\end{remark}

\section{Proof of Theorem~\ref{thm complements line graphs}} \label{sec 5}
Recall that in Theorem~\ref{thm complements line graphs} $G$ is a graph of maximum degree $\Delta\ge 3$, and line graph $H = L(G)$.  Denote for brevity $c(i) = c(i,H)$.

\begin{proof}[Proof of Theorem~\ref{thm complements line graphs}]
The upper bound is a consequence of Theorem~\ref{main thm} and the observation that the complement of the line graph of $G$ is an induced subgraph of $Kn(2,n)$.

In the remainder of the proof, we show that
\begin{equation}\label{eq:central}
    \sum_{i=1}^{|V(H)|-1} c(i) \leq 12\cdot \mathds{1}_{\Delta=3}+16\cdot \mathds{1}_{\Delta=4}+\dfrac{\Delta(\Delta+3)}{2} \mathds{1}_{\Delta\ge 5}.
\end{equation}
The lower bound then directly follows by Lemma~\ref{ACS lemma} applied for $H^c$.

First, since $H$ has maximum degree at most $2(\Delta - 1)$ we have $c(1) \leq 2(\Delta - 1)$. Then, we prove that $c(2) \leq \max(\Delta - 1, 4)$. Indeed, let $\{i_1,i_2\}$ and $\{j_1, j_2\}$ be two distinct vertices of $H$. If the two sets have a common element, assume that $i_1 = j_1$. Then, the common neighbours of $\{i_1,i_2\}$ and $\{j_1, j_2\}$ belong to the set $\{ \{i_1,k\} : v_{i_1} v_k \in E(G) \} \cup \{i_2,j_2\}$. Otherwise, $\{i_1,i_2\}$ and $\{j_1, j_2\}$ have no common element, in which case they can have at most four neighbours:  $\{i_1, j_1\}$, $\{i_1, j_2\}$, $\{i_2, j_1\}$ and $\{i_2, j_2\}$.

We divide the remainder of the proof in two cases according to the value of $\Delta$. Assume first that $\Delta \geq 5$.

\begin{claim}
For every $i$ between $3$ and $\Delta-1$, we have that $c(i)\le \max(\Delta-i, 2)$.
\end{claim}
\begin{proof}
Consider a set $S$ of $i\in [3, \Delta - 1]$ vertices of $H$. If all vertices in $S$ correspond to edges of $G$, containing a fixed vertex, then there are at most $\Delta-i$ vertices of $H$, connected to every vertex in $S$. 

If $S$ contains two vertices of $H$, which correspond to disjoint edges of $G$, say $\{i_1,i_2\}$ and $\{j_1,j_2\}$, there are several cases to consider:
\begin{itemize}
    \item If some of the remaining sets is disjoint from both $\{i_1,i_2\}$ and $\{j_1,j_2\}$, then the vertices in $S$ have no common neighbour in $H$;
    \item If some of the remaining sets intersects both $\{i_1,i_2\}$ and $\{j_1,j_2\}$, then let without loss of generality this set be $\{i_1,j_1\}$. In this case, the vertices in $S$ have at most two common neighbours, corresponding to $\{i_1,j_2\}$ and $\{i_2,j_1\}$;
    \item It remains the case when none of the remaining sets intersects either $\{i_1,i_2\}$ or $\{j_1,j_2\}$. In this case, assume without loss of generality that $\{i_1,k\}$ is in $S$ for some $k \notin \{i_2,j_1,j_2\}$. It means that the vertices in $S$ must have at most two common neighbours in $H$: $\{i_1,j_1\}$ and $\{i_1,j_2\}$.
\end{itemize}
It remains to consider the possibility that all vertices in $S$ correspond to sets, which intersect non-trivially, but do not have a common element. Then one must have $i=3$ and three sets $\{i_1,i_2\}, \{i_2,i_3\}$ and $\{i_1,i_3\}$. In this case the vertices in $S$ have no common neighbour in $H$. This finishes the proof of the claim.
\end{proof}
It remains to note that $c(\Delta) < 2$ since $c(2) \leq \Delta - 1$. Thus, choosing among the neighbours of any given vertex in $H$, one may observe that $c(\Delta) = \dots = c\left( \min(|V(H)|-1, 2(\Delta -1)\right) \leq 1$ and $c(i) = 0$ for every $i > 2(\Delta -1) $. Summing up, we get
\begin{align*}
& \sum_{i=1}^{|V(H)|-1} c(i)\\ 
\leq \hspace{0.3em} 
& 2(\Delta - 1) + \Delta - 1 \ + \sum_{i=3}^{\Delta - 2} (\Delta - i) + 2 + \sum_{i=\Delta}^{2(\Delta -1 )} 1\\ 
=\hspace{0.3em} 
& 3\Delta-3 + \left(\dfrac{(\Delta-3)(\Delta -2)}{2} - 1\right) + 2 + \Delta - 1 \\ 
=\hspace{0.3em}
& \dfrac{\Delta(\Delta + 3)}{2}.
\end{align*}
The first case is proved.

Taking advantage of the first part of the proof, we can mimic the same arguments to obtain upper bounds for the case $\Delta \in \{3, 4\}$. In particular,
\begin{equation*}
    \sum_{i=1}^{|V(H)|-1} c(i) \leq 4+4+2+2 = 12, \ \text{if} \ \Delta=3,
\end{equation*}
and
\begin{equation*}
    \sum_{i=1}^{|V(H)|-1} c(i) \leq 6+4+2+2+1+1 = 16, \ \text{if} \ \Delta=4.
\end{equation*}
Let us explain in more detail the last two cases. We already know that $c(1)\le 2(\Delta-1)$ and $c(2)\le 4$, so by Observation~\ref{ob trivial} $c(i_1) = 0$ for every $i_1\ge 2\Delta - 1$ and $c(i_2)\le 1$ for every $i_2\ge 5$. Thus, it remains to show that $c(3)\le 2$ in both cases, having $c(4)\le c(3)$ then it proves our claim. 

Indeed, on the one hand, if three edges of $G$ have a common endvertex, then at most one other edge may be adjacent to all three since $\Delta\le 4$. On the other hand, if three edges of $G$ are disjoint or if they form a $K_3$, then no edge is adjacent to all of them. It remains the case when two of the edges (say $uv, uw\in E(G)$) have a common endvertex (in this case $u$), not adjacent to the third edge. Then, if $\Delta = 3$, the only edges of $G$ that may possibly be adjacent to all three edges are $vw$ and the third edge of $G$, containing $u$. If $\Delta = 4$, there may possibly be at most three edges of $G$ of the form $uv_1, uw_1$ and $vw$, adjacent to both $uv$ and $uw$, but no edge of $G$ different from $uv$ and $uw$ may be adjacent to $uv_1, uw_1$ and $vw$ at the same time.

This concludes the proof of the theorem.
\end{proof}


\section{Conclusion and further questions}\label{sec: conclusion}
In this paper we studied the boxicity of Kneser graphs. Finding the right value of the boxicity of $Kn(k,n)$ for every choice of $k,n$ with $n\ge 2k+1$ seems an interesting, but also quite challenging problem. An easier, but nonetheless intriguing question is whether the lower bound on $Kn(k,n)$ in Theorem \ref{main thm} could be improved to $n - \Theta(1)$ when $n\geq C k^3$ for a large enough constant $C$.
Moreover, in the case $k=2, n\ge 5$, we proved that $\boxi(Kn(2,n))$ is either $n-3$ or $n-2$, and with the help of a \textit{SageMath} program we could show that the boxicty of the Petersen graph, corresponding to $Kn(2,5)$, is 3. This result suggests that the right value of $\boxi(Kn(2,n))$ might be $n-2$ for every $n\ge 5$.

\section{Acknowledgements}
The authors would like to thank Mat\v{e}j Stehl\'ik for a number of useful discussions.

\bibliographystyle{plain}
\bibliography{References}

\end{document}